\documentclass[journal,10pt]{IEEEtran}
\usepackage{ifpdf}

\ifCLASSINFOpdf
\else
\fi
\usepackage{amsmath}
\usepackage{algorithmicx}
\usepackage{algpseudocode}

\usepackage{array}

\usepackage{amsfonts}
\usepackage{amssymb}
\usepackage{graphicx}
\usepackage{float}
\usepackage{setspace}

\begin{document}
%
\title{Maximum likelihood estimation of covariances of elliptically symmetric distributions}


\author{\IEEEauthorblockN{Christophe Culan 
		,
Claude Adnet\IEEEauthorrefmark{1}, 
}
\IEEEauthorblockA{ 
	Advanced Radar Concepts division,
Thales, Limours 91338, France}
\thanks{Manuscript received October XX, 2016; revised October XX, 2016. 
Corresponding author: M. Culan (email: christophe.culan@thalesgroup.com).}}

\markboth{ArXiv prepublication}
{Culan \MakeLowercase{\textit{et al.}}: Maximum likelihood estimation of covariances of elliptically symmetric distributions}
%



\IEEEtitleabstractindextext{%
\begin{abstract}
	Elliptically symmetric distributions are widely used in portfolio modeling, as well as in signal processing applications for modeling impulsive background noises. Of particular interest are algorithms for covariance estimation and subspace detection in such backgrounds.\\
	This article tackles the issue of correctly estimating the covariance matrix associated to such models and detecting additional signal superimposed on such distributions. A particular attention is given to the proper accounting of the circular symmetry for the subclass of complex elliptical distributions in the case of complex signals.\\
	In particular Tyler's estimator is shown to be a maximum likelihood estimate over all elliptical models, and its extension to the complex case is shown to be a maximum likelihood estimate for the subclass of complex elliptical models (CES); other M-estimators are also shown to be maximum likelihood estimates over some restricted classes of elliptical models. The extension of Tyler's and other M-estimators to constrained covariance estimation is also discussed, in particular for toeplitz constrains.\\
	Finally likelihood ratio signal detection tests associated to the various estimators introduced in this article are also discussed.
\end{abstract}

}

\maketitle

\IEEEdisplaynontitleabstractindextext

%
\IEEEpeerreviewmaketitle

\section{Introduction: elliptically symmetric distributions}
\label{sec:intro}

Elliptically symmetric distributions constitute a wide class of distributions which generalize multivariate gaussian models \cite{chmielewski1981elliptically}\cite{van1995multivariate}\cite{ollila2012complex}. Some examples are the K-distributions \cite{conte1991modelling}, the t-distributions \cite{krishnaiah1986complex} \cite{ollila2003robust}, the generalized gaussian distributions \cite{novey2010complex}, or the larger class of compound gaussian models (CG) \cite{conte1987characterisation} \cite{gini2002vector} \cite{pascal2008covariance}, which are widely used in simulations as they are easy to generate.\\
Such distribution models are generally applied for applications in which the symmetry properties of a standard multivariate gaussian distribution are desirable, but for which one wants to modelize heavy tailed distributions. Some applications of interest include portfolio modeling in statistical finance, as well as the modeling of highly impulsive noises in signal processing such as clutter signal in RADAR and SONAR applications\cite{conte1995asymptotically}.

\subsection{Notations and conventions, and some useful formulae}
\label{ssec:notations}

In the following development, the following notations and conventions shall be observed:
\begin{itemize}
	\item For any topological space $\mathcal{R}_1$ and $\mathcal{R}_2$, the topological space $\mathcal{R}_1 \times \mathcal{R}_2$ is the product of $R_1$ and $R_2$, whereas $\mathcal{R}_1 \vee \mathcal{R}_2$ is the disjoint union of $R_1$ and $R_2$. One also notes, for a topological space $\mathcal{R}$,  $\mathcal{R}^N$ to be the product of $N$ copies of $\mathcal{R}$ and $\vee^N \mathcal{R}$ to be the disjoint union of $N$ copies of $\mathcal{R}$.
	\item Let $\mathcal{R}_1$, ..., $\mathcal{R}_N$ be topological spaces, and let $\mu_n$ be a measure defined on $\mathcal{R}_n$ for $1 \leq n \leq N$. The measure $\mu_1 \dots \mu_n = \prod_{n=1}^N \mu_n$ is the product measure of $(\mu_n)_{1\leq n\leq N}$ defined on $\prod_{n=1}^N \mathcal{R}_n$, whereas\\ $\mu_1 \vee \dots \vee \mu_n = \bigvee_{n=1}^N \mu_n$ is the measure of $\bigvee_{n=1}^N \mathcal{R}_n$ such that its restriction to $\mathcal{R}_n$ is $\mu_n$ for $1 \leq n \leq N$ and is the joint measure of $(\mu_n)_{1 \leq n \leq N}$.\\
	Moreover for a topological space $\mathcal{R}$ and a measure $\mu$ defined on $\mathcal{R}$, $\mu^N$ is the product measure of $N$ copies of $\mu$ defined on $\mathcal{R}^N$, whereas $\vee^N \mu$ is the joint measure of $N$ copies of $\mu$ defined on $N\cdot\mathcal{R}$.
	\item The vector space $\mathbb{C}^d$ is canonically identified to $\mathbb{R}^{2d}$ \cite{scharf1991statistical}.
	\item $\mathcal{H}_d(\mathbb{R})$ is the set of symmetric matrices of size $(d,d)$; $\mathcal{H}_d(\mathbb{C})$ is the set of hermitian matrices of size $(d,d)$.\\
	$\mathcal{H}_d^+(\mathbb{K})$ is the subset of matrices of 	$\mathcal{H}_d^+(\mathbb{K})$ which are positive; $\mathcal{HP}_d^+(\mathbb{K})$ is the subset of matrices of $\mathcal{H}_d^+(\mathbb{K})$ of unit determinant.
	\item  $X^\dagger$ is the transpose of any real matrix or vector $X$, or the transpose conjugate of any complex matrix or vector $X$.
	\item $\mu_{\mathbb{R}^d}$ is the Lebesgues measure of the space $\mathbb{R}^d$. Similarly, $\mu_{\mathbb{C}^d}$ is the Lebesgues measure of $\mathbb{C}^d \simeq \mathbb{R}^{2d}$.
	\item $\mu_{\mathbb{R}_+}$ designates the restriction of $\mu_\mathbb{R}$ to $\mathbb{R}_+$.
	\item $\mathcal{S}_{d-1}$ is the $(d-1)$-sphere, which is identified as the following part of $\mathbb{R}^d$: $\mathcal{S}_{d-1} \simeq \left\{ x \in \mathbb{R}^d; x^\dagger x = 1 \right\}$. 
	\item $s_{d-1}$ shall denote the probability distribution of $\mathcal{S}_{d-1}$ isotropic for the canonical scalar product of $\mathbb{R}^d$, defined by:
	\begin{equation}
	\left<y,x\right> = \sum\limits_{k=0}^{d-1} y_k x_k
	\end{equation}
	\item $\mathcal{U}_{d}$ designates $S_{2d-1}$ identified to $\left\{ x \in \mathbb{C}^d; x^\dagger x = 1 \right\}$ in $\mathbb{C}^d$.
	\item $u_{d}$ designates the probability measure of $\mathcal{U}_d$ isotropic for the canonical scalar product of $\mathbb{C}^d$, defined by:
	\begin{equation}
	\left<y,x\right> = \sum\limits_{k=0}^{d-1} \overline{y_k}x_k
	\end{equation}
	\item $[\mathbb{C}^d]$ designates the space $\frac{\mathbb{C}^d}{\mathcal{U}_1}$, whose elements are noted:\\  $[x] = \left\{ y\in\mathbb{C}^d; \exists u \in \mathcal{U}_1, y = ux \right\}$. It is a space of dimension $2d-1$ over $\mathbb{R}$; one defines its Lebesgues measure $\mu_{[C^d]}$ to be its unique measure $\mu$ invariant by translation such that $\int_{[x]\in S_{[\mathbb{C}^d]}^\infty} d\mu\left([x]\right) = 1$ with $S_{\mathbb{C}^d}^\infty = \left\{[x] \in [\mathbb{C^d}]; \max_{0 \leq k \leq d-1} \left| x_k\right|\leq 1\right\}$.
	\item $\mathbb{C}P^{d-1}$ designates the $d-1$ projective space over $\mathbb{C}$; Note that  $\mathbb{C}P^{d-1} \simeq \frac{\mathcal{U}_d}{\mathcal{U}_1}$; its elements are noted:\\ $[x] = \left\{u x; u \in \mathcal{U}_1 \right\}$ with $x \in \mathcal{U}_d$.
	\item $\delta_x^\mathcal{R}$ is the Dirac delta distribution centered in $x$ in the space $\mathcal{R}$.
	\item Measures on $\mathbb{R}^d$ or $\mathbb{C}^d$ and their respective spheres can be related in the following ways:
	\begin{equation}
	\label{eq:meas_lebegues_sphere}
	\begin{array}{l}
	d\mu_{\mathbb{R}^d}(x) = \frac{2\pi^{\frac{d}{2}}}{\Gamma\left(\frac{d}{2}\right)}{\left\|x\right\|}^{d-1} d\mu_{\mathbb{R_+}}\left(\left\|x\right\|\right) s_{d-1}\left(\frac{x}{\left\|x\right\|}\right)\vspace{1.5 mm}\\
	d\mu_{\mathbb{C}^d}(x) = \frac{2\pi^d}{(d-1)!}{\left\|x\right\|}^{2d-1} d\mu_{\mathbb{R_+}}\left(\left\|x\right\|\right) u_d\left(\frac{x}{\left\|x\right\|}\right)
	\end{array}
	\end{equation}
	\item One shall note:
	\begin{equation}
	\label{eq:def_entropy}
	h_\delta = H\left(\delta_0^{\mathbb{R}}|\mu_{\mathbb{R}}\right) = \int_{t \in \mathbb{R}} \log\left(\frac{d\delta_0^{\mathbb{R}}}{dt}(t)\right)d\delta_0^{\mathbb{R}}(t)
	\end{equation}
	$h_\delta$ is the entropy of the Dirac distribution relative to the Lebesgues measure of $\mathbb{R}$, which is considered as positive and infinite. It can be used to express different other entropies in higher dimensional settings:\\
	\begin{equation}
	\begin{array}{l}
	\label{eq:entropy_higher_d}
	H\left(\delta_x^{\mathbb{R}^d}|\mu_{\mathbb{R}^d}\right) = d h_\delta\vspace{1.5 mm}\\
	H\left(\delta_x^{\mathbb{C}^d}|\mu_{\mathbb{C}^d}\right) = 2d h_\delta
	\end{array}
	\end{equation}
\end{itemize}

\subsection{A first step: Spherically symmetric distributions as a generalization of isotropic gaussian distributions}
\label{ssec:spherical}

A spherically symmetric distribution of $\mathbb{R}^d$ is a distribution which, like the isotropic gaussian distribution, exhibits complete invariance by any linear orthonormal transformation. Such a distribution $P$ can be noted:

\begin{equation}
\label{eq:def_spherical}
dP(x) = dQ\left(\left\| x \right\|\right)ds_{d-1}\left(\frac{x}{\left\| x \right\|}\right)
\end{equation}

with $Q$ being any probability measure on $\mathbb{R}_+$ \cite{ollila2012complex}. $Q$ is the radial distribution of the spherically symmetric distribution $P$.\\
Note for example that an isotropic gaussian distribution would correspond to $dQ(r) = \sigma dQ_0(\sigma r)$ with $Q_0$ being a $\chi(d)$ distribution and $\sigma \in \mathbb{R}_+$ \cite{scharf1991statistical}.\\
Since the distribution $Q$ can take on any value, this allows for the modeling of arbitrarily impulsive isotropic distributions.\\

\subsection{Application to the anisotropic case : Elliptically symmetric distributions}
\label{ssec:elliptical}
Spherically symmetric distributions can be generalized to take into account distributions which would be isotropic for any scalar product over $\mathbb{R}^d$ by applying a simple linear variable change:

Let $A \in GL_d\left(\mathbb{R}\right)$ and $X$ be a random variable following a spherically symmetric distribution of radial distribution $Q$. Then the random variable $A X$ follows a full rank elliptically symmetric distribution \cite{chmielewski1981elliptically}\cite{ollila2012complex}. Its distribution can be recast as:

\begin{equation}
\label{eq:def_elliptical}
dP(x) = \left|A\right|dQ\left(\left\|A x\right\|\right)ds_{d-1} \left( \frac{A x}{\left\| A x \right\|} \right)
\end{equation}

Noting $R = \left(A^\dagger A\right)^{-1}$ which is called the correlation matrix of the distribution $P$ and noting its Cholesky decomposition $R = L(R)L(R)^\dagger$, this distribution can be recast as:

\begin{equation}
\label{eq:def_elliptical_RQ}
dP(x) = \frac{1}{\sqrt{\left|R\right|}}dQ\left(\sqrt{x^\dagger R^{-1}x}\right)ds_{d-1}\left( \frac{L(R)^{-1} x}{\sqrt{x^\dagger R^{-1} x}} \right)
\end{equation}

Thus an elliptically symmetric distribution (ES) is completely parametrized by its radial distribution $Q$ and its correlation matrix $R$. Moreover those are uniquely defined if one adds a normalizing constrain to $R$, such as $|R|=1$ or $\text{tr}(R)=\text{Cste}$. The condition $|R|=1$  is generally assumed in theoretical studies for convenience; indeed the expression as given by \eqref{eq:def_elliptical_RQ} then simplifies as:

\begin{equation}
\label{eq:def_elliptical_normRQ}
dP(x) = dQ\left(\sqrt{x^\dagger R^{-1}x}\right)ds_{d-1}\left( \frac{L(R)^{-1} x}{\sqrt{x^\dagger R^{-1} x}} \right)
\end{equation}

On the other hand the condition $\text{tr}(R)=\text{Cste}$ is much easier to handle in numerical computations.\\

The class of ES distributions is quite vast and in many applications, one could restrict the study to elliptically symmetric distributions of a specific kind, for which the radial distribution is parametrized by a few parameters. In this paper, we focus on two following situations:

\begin{itemize}
	\item  the case when one has no a priori on the radial distribution. The radial distribution can therefore be any probability distribution over $\mathbb{R}_+$. This is the most general case.
	\item the case in which the radial distribution is a scaled version of a unique base distribution: $dQ(r) = \sigma dQ_0(\sigma r)$, with $\sigma >0$; this is notably the case of gaussian distributions, for which $Q$ corresponds to a $\chi^2(d)$ distribution. Thus such a class of distributions is parametrized by its correlation matrix $R$ and the positive scale factor $\sigma$.
\end{itemize}

\section{Covariance estimation of an ES distribution}
\label{sec:cov_est}

This section deals with the derivation of estimators of the covariance matrix for elliptically symmetric distributions.

\subsection{Maximum likelihood: general motivation and background}
\label{ssec:ml}

The principle of maximum likelihood estimation is applied in many fields, as it generally has many desirable properties \cite{pfanzagl1994parametric}. As was noted by Akaike \cite{akaike1998information}, this principle can in fact be restated in terms of entropy.\\

Indeed consider a probability distribution model $P$ and a sampling distribution $S$ defined on some space $\mathcal{R}$, which corresponds to the observed data. The log-likelihood score of the model $P$ for the distribution model $S$ can be defined by:

\begin{equation}
\label{eq:logl_def}
l(P|S)= \int_{x \in \mathcal{R}} \log\left(\frac{dP}{dS}(x)\right)dS(x)
\end{equation}

This in fact corresponds to minus the entropy of the sampling distribution S relative to the model distribution P (also known as KL-divergence) \cite{kullback1951information}:

\begin{equation}
\label{eq:logl_entropy}
l(P|S) = -H(S|P)	
\end{equation}

Hence maximizing the likelihood $l(P|S)$ is equivalent to minimizing the amount of information lost by choosing the model $P$ compared to the information contained in the sampling distribution $S$.

\subsection{Maximum likelihood under independent, identically distributed hypothesis}
\label{ssec:ml_iid}

A very common hypothesis in estimation problems is that of supposing that one has $N$ samples which are a priori identically distributed, and independant of each other.\\

The considered distributions are then valued in $\mathcal{R} = \vee^N \mathcal{B}$, with $\mathcal{B}$ being some base space; the model distribution is of the form $\frac{1}{N}\vee^N\mathcal{P}$ whereas the sampling distribution can be noted $S = \frac{1}{N} \bigvee_{n=1}^N S_n$, with $S_n$ being defined on $\mathcal{B}$ for $1 \leq n \leq N$ and corresponding to a single sample.

The likelihood can therefore be expressed as:

\begin{equation}
\label{eq:logl_mean}
l(P|S) = \frac{1}{N}\sum\limits_{n=1}^N l(P|S_n)
\end{equation}

This means that the likelihood for $N$ i.i.d samples is the average of the likelihood for each sample taken separately.\\

Now under no specific hypothesis other than the fact that each sample $x_n$ is in $\mathcal{B}$, the sampling distribution is given by: $S = \bigvee_{n=1}^N \delta_{x_n}^\mathcal{B}$ and the likelihood is given by:

\begin{equation}
\label{eq:logl_mean_delta}
\begin{split}
l(P|S) = &\frac{1}{N}\bigvee_{n=1}^N l\left(P|\delta_{x_n}^\mathcal{B}\right)\\ =& \frac{1}{N}\sum\limits_{n=1}^N\log\left(\frac{dP}{d\mu}(x_n)\right) -\frac{1}{N}\sum\limits_{n=1}^N H\left(\delta_{x_n}^\mathcal{B}|\mu\right)
\end{split}
\end{equation}

with $\mu$ being any volume distribution of $\mathcal{B}$. Thus one finds back the usual expression of the likelihood up to some constant terms.

\subsubsection{Likelihood of an ES distribution}
\label{sssec:logl_es}

Let us now consider an ES distribution $P$ defined on $\mathbb{R}^d$, parametrized by its radial distribution $Q$ and its correlation matrix $R$ constrained to $|R| = 1$, and a sampling distribution $S$ on $\mathbb{R}^d$.\\

$S$ can be decomposed in $R$-elliptical coordinates: 
\begin{equation}
\label{eq:sampling_elliptical}
dS(x) = dQ_S\left(\sqrt{x^\dagger R^{-1} x}\right) d\widehat{S}_{\sqrt{x^\dagger R^{-1} x}}\left(\frac{L(R)^{-1} x}{\sqrt{x^\dagger R^{-1} x}}\right)
\end{equation}

with $Q_S$ being a probability distribution over $\mathbb{R}_+$, and $\widehat{S}_r$ being a probability distribution of the unit sphere $\mathcal{S}_{d-1}$ for any $r \in \mathbb{R}_+$.\\

The log-likelihood of the distribution $P$ for a sampling distribution $S$ can then be recast by separating the integration over $r = \sqrt{x^\dagger R^{-1} x} \in \mathbb{R}_+$ and $u = \frac{L(R)^{-1} x}{\sqrt{x^\dagger R^{-1} x}} \in \mathcal{S}_{d-1}$:

\begin{flalign*}
l(Q,R|S) &&
\end{flalign*}
\begin{equation}
\label{eq:logl_elliptical}
\begin{split}
=&\iint\limits_{
	\substack{
	r \in \mathbb{R}_+\\
	u \in \mathcal{S}_{d-1}
	}
	}
 \!\!\! \log\!\left(\frac{1}{r^{d-1}}\frac{dQ}{dQ_S}(r)r^{d-1} \frac{ds_{d-1}}{d\widehat{S}_r}(u)\right) \! d\widehat{S}_r(u)  dQ_S(r)\\\\
 =&\iint\limits_{
 	\substack{
 		r \in \mathbb{R}_+\\
 		u\in \mathcal{S}_{d-1}
 	}
 } \log\left(r^{d-1}\frac{ds_{d-1}}{d\widehat{S}_r}(u)  \right) d\widehat{S}_r(u)dQ_S(r)\\
 &+\int_{r \in \mathbb{R}_+} \log\left(\frac{dQ}{dQ_S}(r)\right) dQ_S(r)\\
 &-(d-1) \int_{r \in \mathbb{R}_+} \log(r) dQ_S(r)\\\\
 =&\iint\limits_{
 	\substack{
 		r \in \mathbb{R}_+\\
 		u\in \mathcal{S}_{d-\!1}
 	}
 } \log\left(r^{d-1}\frac{ds_{d-1}}{d\widehat{S}_r}(u)\right) d\widehat{S}_r(u) dQ_S(r)\\
 &- H(Q_S|Q)-(d-1) \int_{r \in \mathbb{R}_+} \log(r) dQ_S(r)
\end{split}
\end{equation}

In the case of a standard sampling distribution corresponding to a single sample $x_0$ in $\mathbb{R}^d$, the corresponding sampling distribution can be expressed in $R$-elliptical coordinates as:

\begin{equation}
\label{eq:delta_elliptical}
d\delta_{x_0}^{\mathbb{R}^d} \! = dQ_{x_0}\left(\sqrt{x^\dagger R^{-1} x}\right)d\widehat{S_{x_0}}_{\sqrt{x^\dagger R^{-1} x}}\left(\frac{L(R)^{-1} x}{\sqrt{x^\dagger R^{-1} x}}\right)
\end{equation}

with:

\begin{equation}
\label{eq:delta_split}
\left\{ \begin{array}{l}
Q_{x_0} =  \delta_{\sqrt{{x_0}^\dagger R^{-1} x_0}}^{\mathbb{R}_+}\\
\widehat{S_{x_0}}_r = \delta_{\frac{L(R)^{-1} x_0}{\sqrt{{x_0}^\dagger R^{-1} x_0}}}^{\mathcal{S}_{d-1}}
\end{array} \right.
\end{equation}

Hence the likelihood of a single sample can be obtained from equation (\ref{eq:logl_elliptical}):

\begin{equation}
\label{eq:logl_delta_elliptical}
\begin{split}
l\left(Q,R|\delta_{x_0}^{\mathbb{R}^d}\right) =  -&(d-1)h_\delta+\log\left(\frac{\Gamma(d)}{2\pi^\frac{d}{2}}\right)\\
-&  H\left(\delta_{\sqrt{{x_0}^\dagger R^{-1} x_0}}^{\mathbb{R}_+}|Q\right)\\
-&\frac{d-1}{2} \log\left({x_0}^\dagger R^{-1} x_0\right)
\end{split}
\end{equation}

Thus for $N$ i.i.d samples $\left(x_n\right)_{1 \leq n \leq N}$, the corresponding likelihood is given by:

\begin{equation}
\label{eq:logl_delta_elliptical_mean}
\begin{split}
l\left(Q,R|\frac{1}{N}\bigvee_{n=1}^N \delta_{x_n}^{\mathbb{R}^d}\right) = &
-(d-1)h_\delta+\log\left(\frac{\Gamma(d)}{2\pi^\frac{d}{2}}\right)\\
&  -\frac{1}{N}\sum\limits_{n=1}^{N} H\left(\delta_{\sqrt{{x_n}^\dagger R^{-1} x_n}}^{\mathbb{R}_+}|Q\right)\\
&-\frac{d-1}{2N} \log\left({x_n}^\dagger R^{-1} x_n\right)
\end{split}
\end{equation}

This can be first maximized with respects to $Q$, leading to:

\begin{equation}
\label{eq:elliptical_ml_Q}
Q = \frac{1}{N} \sum\limits_{n=1}^N \delta_{\sqrt{x_n R^{-1} x_n}}^\mathbb{R_+}
\end{equation}

This results in the following concentrated likelihood, maximized over all possible radial distributions for a given correlation matrix:

\begin{equation}
\label{eq:elliptical_concentrated}
\begin{split}
l\left(R|\frac{1}{N}\bigvee_{n=1}^N \delta_{x_n}^{\mathbb{R}^d}\right) =& -(d-1)h_\delta-\log\left(\frac{\Gamma(d)}{2\pi^\frac{d}{2}}\right)\\
&-\frac{d-1}{2N} \sum\limits_{n=1}^N \log\left({x_n}^\dagger R^{-1} x_n\right)
\end{split}
\end{equation}

\subsubsection{Likelihood of a CES distribution}
\label{sssec:logl_ces}

Similarly to the case of an elliptically symmetric distribution, a complex elliptically symmetric distribution can be expressed as:

\begin{equation}
\label{eq:def_ces}
dP(x) = dQ\left(\sqrt{x^\dagger R^{-1} x}\right)du_{d}\left(\frac{L(R)^{-1}x}{\sqrt{x^\dagger R^{-1}x}}\right)
\end{equation}

with $R$ being a hermitan, positive definite, unit determinant matrix and $Q$ verifying that $Q\left(\mathbb{R}_+\right)=1$.

Thus the likelihood of such a model $P$ for any given sampling distribution $S$ can be expressed as:

\begin{samepage}
\begin{flalign*}
l(Q,R|S) &&
\end{flalign*}
\begin{equation}
\label{eq:logl_ces}
\begin{split}
 = & \iint\limits_{
 	\substack{
		r \in \mathbb{R}_+\\
		u \in \mathcal{U}_d
		}
		}
		\log\left(\frac{1}{r^{2(d-1)}}\frac{dQ}{dQ_S}(r)r^{2(d-1)} \frac{du_d}{d\widehat{S}_r}(u)\right) d\widehat{S}_r(u)dQ_S(r)\\\\
=&\iint\limits_{\substack{
	r \in \mathbb{R}_+\\
	u\in \mathcal{U}_d
	}
	} \log\left(r^{2(d-1)}\frac{du_d}{d\widehat{S}_r}(u)\right) d\widehat{S}_r(u)dQ_S(r)\\
&+\int_{r \in \mathbb{R}_+} \log\left(\frac{dQ}{dQ_S}(r)\right) dQ_S(r)\\
&- 2(d-1) \int_{r \in \mathbb{R}_+} \log(r) dQ_S(r)\\\\
=&\iint\limits_{\substack{	
	r \in \mathbb{R}_+\\
	u\in \mathcal{U}_d
	}
}\log\left(r^{2(d-1)}\frac{du_d}{d\widehat{S}_r}(u)\right) d\widehat{S}_r(u)dQ_S(r)\\
&- H(Q_S|Q)-2(d-1) \int_{r \in \mathbb{R}_+} \log(r) dQ_S(r)
\end{split}
\end{equation}
\end{samepage}

with $S$ admitting the following decomposition in $R$-elliptical coordinates as:

\begin{equation}
\label{eq:sampling_ces}
dS(x) = dQ_S\left(\sqrt{x^\dagger R^{-1} x}\right) d\widehat{S}_{\sqrt{x^\dagger R^{-1} x}}\left(\frac{L(R)^{-1} x}{\sqrt{x^\dagger R^{-1} x}}\right)
\end{equation}

Now let us consider a single sample $x_0 \in \mathbb{C^d}$. If no other hypothesis is done on $x_0$, the corresponding sampling distribution would be given by $\delta_{x_0}^{\mathbb{C}^d}$; however an implicit hypothesis when one works in $\mathbb{C}^d$ rather than $\mathbb{R}^{2d}$ is that the signal exhibits circular symmetry, meaning that it is phase invariant. The sampling distribution associated to a single sample under this hypothesis is given by:

\begin{equation}
\label{eq:def_csampling}
\delta_{[x_0]}^{\mathbb{C}^d} = u_1\delta_{[x_0]}^{[\mathbb{C}^d]}
\end{equation}

This sampling distribution can be decomposed in $R$-elliptical coordinates as:

\begin{equation}
\label{eq:delta_ces}
d\delta_{[x_0]}^{\mathbb{C}^d} \! = \! dQ_{[x_0]} \!\! \left(\sqrt{x^\dagger R^{-1} x}\right) \! d\widehat{S_{[x_0]}}_{\sqrt{x^\dagger R^{-1} x}} \! \left(\frac{L(R)^{-1} x}{\sqrt{x^\dagger R^{-1} x}}\right)
\end{equation}

with:

\begin{equation}
\label{eq:delta_csplit}
\left\{ \begin{array}{l}
Q_{[x_0]}(r) = \delta_{\sqrt{{x_0}^\dagger R^{-1} x_0}}^{\mathbb{R}_+}\\
\widehat{S_{[x_0]}}_r =
u_1\delta_{\left[\frac{L(R)^{-1}x}{\sqrt{x^\dagger R^{-1} x}}\right]}^{[\mathcal{U}_{d}]} = \delta_{\left[\frac{L(R)^{-1}x}{\sqrt{x^\dagger R^{-1} x}}\right]}^{\mathcal{U}_{d}}
\end{array} \right.
\end{equation}

Thus the likelihood with respect to such a distribution for $N$ i.i.d samples $\left(x_n\right)_{1 \leq n \leq N}$ can be expressed as:

\begin{equation}
\label{eq:logl_delta_ces}
\begin{split}
l\left(Q,R|\bigvee_{n=1}^N \delta_{[x_n]}^{\mathbb{C}^d}\right) = & -2(d-1)h_\delta+\log\left(\frac{(d-1)!}{\pi^{d-1}}\right)\\
&  -\frac{1}{N} \sum\limits_{n=1}^N H\left(\delta_{\sqrt{{x_n}^\dagger R^{-1} x_n}}^{\mathbb{R}_+}|Q\right)\\
&-\frac{(d-1)}{N} \sum\limits_{n=1}^N \log\left({x_n}^\dagger R^{-1} x_n\right)
\end{split}
\end{equation}

By maximizing this likelihood with respect to the radial distribution, one gets:

\begin{equation}
\label{eq:ces_concentrated}
\begin{split}
l(R|S) = & -2(d-1)h_\delta+\log\left(\frac{(d-1)!}{\pi^{d-1}}\right)\\
&-\frac{(d-1)}{N} \sum\limits_{n=1}^N \log\left({x_n}^\dagger R^{-1} x_n\right)
\end{split}
\end{equation}	

Hence this is the same expression as in the real case, up to a multiplicative and additive constant. Ignoring the additive constants, the likelihood corresponds to the following function:

\begin{equation}
\label{eq:general_es_concentrated}
l(R|S) = -c_\mathbb{K}\frac{d-1}{2N}\sum\limits_{n=1}^N \log\left({x_n}^\dagger R^{-1} x_n\right)
\end{equation}

with $c_\mathbb{K}=1$ for $\mathbb{K} = \mathbb{R}$, and $c_\mathbb{K}=2$ for $\mathbb{K} = \mathbb{C}$.

\subsubsection{Maximum likelihood in the unconstrained case: Tyler's estimator}
\label{subsub:tyler}

As is shown in \cite{wiesel2012geodesic}, the resulting concentrated likelihood is convex on the set $\mathcal{HP}_d^+(\mathbb{K})$, equipped with the metric: $d(R_1,R_2) = \text{tr}\left(\log^2(R_1^{-1}R_2)\right)$; hence this guaranties the existence of a unique maximizer of this function.\\

Moreover one can extract an implicit equation which defines this maximizer. Indeed, by introducing a Lagrange multiplier to ensure the determinant constrain on $R$ and ignoring the constant terms, one gets:

\begin{equation}
\label{eq:logl_lagrange}
l_\lambda(R) = -\lambda \log|R| -c_\mathbb{K}\frac{d-1}{2N} \sum\limits_{n=1}^N \log\left({x_n}^\dagger R^{-1} x_n\right)
\end{equation}

This function can be differentiated in $R$, yielding the following maximum likelihood equation:

\begin{equation}
\label{eq:es_ml_eq}
\text{tr}\left(\left(\lambda R-c_\mathbb{K}\frac{d-1}{2N} \sum\limits_{n=1}^N \frac{x_n {x_n}^\dagger}{{x_n}^\dagger R^{-1} x_n}\right)d\left(R^{-1}\right)\right)=0
\end{equation}

With no other constrains on the correlation matrix, this results in the following implicit equation, which defines $R$ up to a constant:

\begin{equation}
\label{eq:es_ml_eq_noconstrain}
R = c_\mathbb{K}\frac{d-1}{2\lambda N} \sum\limits_{n=1}^N \frac{x_n {x_n}^\dagger}{{x_n}^\dagger R^{-1} x_n}
\end{equation}

This equation has already been obtained using different model hypothesis, using angular distributions $\cite{auderset2005angular}$ \cite{tyler1987statistical} or a compound gaussian model background $\cite{pascal2008covariance}$. Although it cannot be explicitly solved, the Tyler's fixed point algorithm can be used to approach its unique solution up to any desirable precision \cite{tyler1987statistical}\cite{pascal2008covariance}.\\


\subsubsection{Constrained estimation of the correlation matrix}
\label{sssec:constrained_estimation}

We are now interested in the case in which the correlation matrix is constrained to be on a closed submanifold $\mathcal{M} \subset \mathcal{HP}_d^{+}\left(\mathbb{K}\right)$. Such cases are in general difficult to solve, as the concentrated likelihood is not necessarily convex. However in general, one can produce a normalized version of an estimator $e$ defined for gaussian models under the following constrain on the covariance matrix : $\Sigma \in \left\{\sigma\mathcal{M}; \sigma \in \mathbb{R}_+\right\}$\\

\begin{onehalfspace}
\begin{samepage}
\begin{algorithmic}[1]
	\Function{Tyler\_of}{$e$,${x_n}_{1 \leq n \leq N},\epsilon,K_\text{max}$}
		\State $R_{-1} \gets I$
		\For{$k$ \textbf{from} $1$ \textbf{to} $K_\text{max}$}
			\For{$n$ \textbf{from} $1$ \textbf{to} $N$}
				\State $y_n \leftarrow \frac{x_n}{\sqrt{{x_n}^\dagger R_{-1} x_n}}$\vspace{1.5 mm}
			\EndFor
			\State $S \gets e\left(\left(y_n\right)_{1 \leq n \leq N}\right)$\vspace{1.5 mm}
			\State $R \gets \frac{S}{\text{tr}(S)}$\vspace{1.5 mm}
			\If{$\text{tr}\left({\left(R_{-1}R-I\right)}^{2}\right) \leq \epsilon$}\vspace{1.5 mm}
				\State \textbf{Break}
			\Else
				\State $R_{-1} \leftarrow R^{-1}$
			\EndIf
		\EndFor
		\State \Return $R_{-1}$
	\EndFunction
\end{algorithmic}
\end{samepage}
\end{onehalfspace}
\vspace{2 mm}

The convergence of this procedure is not insured in general; therefore this should be checked depending on the imposed constrain and the estimator $e$ which is used. However if it converges and $e$ is a maximum likelihood estimator of the covariance matrix for centered gaussian models, then the corresponding normalized estimator is a solution to the maximum likelihood equation; it is not guaranted, however, to be the global maximum of the likelihood function.

\subsubsection{Application to the Toeplitz constrain for CES distributions: the Burg-Tyler algorithm}
\label{sssec:toeplitz_estimation}

Based on the procedure described above, we propose an algorithm for estimating the correlation matrix of a stationary, regularly sampled complex signal.\\
Unfortunately the Toeplitz constrain is quite difficult to enforce even in the simpler case of gaussian models; indeed it fails to have a convex structure in this case, and thus the existence of a unique solution to the ML equation is not guarantied. We propose instead to used the previously described procedure on the well known Burg method \cite{burg1978maximum}\cite{ulrych1976time}.\\ 

%

The Burg algorithm estimates the Schur coefficients $\left(\mu_m\right)_{1 \leq m \leq d-1}$ and the residual variance $\sigma^2$ associated with the corresponding AR process. 
\\
These Schur coefficients together with the residual variance $\sigma^2$ can be used to obtain the covariance matrix by using Trench's algorithm \cite{trench1964algorithm}\cite{zohar1969toeplitz}\cite{barbaresco1997analyse}, or by first obtaining the coefficients of the corresponding AR polynomial \cite{decurninge2015quantiles} \cite{lehmer1961machine}, and then by using the Gohberg-Semencul formula \cite{gohberg2005convolution} \cite{kailath1978inverses} \cite{landau1987maximum} \cite{gohberg1972inversion}.\\

%
%
%
The Burg method has since been generalized into a multisegment Burg algorithm, which allow to take into account several independent samples for performing the estimation \cite{haykin1982maximum}. We therefore propose to use the normalization procedure introduced in ~\ref{sssec:constrained_estimation} on this multisegment Burg method, which offers reasonable performances for a limited computation cost; one then obtains the following Burg-Tyler algorithm:\\

%

\begin{onehalfspace}
	\begin{samepage}
		\begin{algorithmic}[1]
			\Function{BT}{$\left(x_n\right)_{1 \leq n \leq N},\epsilon,K_{\max}$}
			\State $\mu \gets (0)_{1 \leq m \leq d-1}$
			\For{$k$ \textbf{from} 1 \textbf{to} $K_{\max}$}
			\State $R_{-1} \gets \textproc{Trench}(1,\mu)$  
			\For{$n$ \textbf{from} 1 \textbf{to} $N$}
			\State $y_n \gets \frac{x_n}{\sqrt{{x_n}^\dagger R_{-1} x_n}}$\vspace{1.5 mm}
			\EndFor\vspace{1.5 mm}
			\State $\left(\sigma^2,\nu\right) \gets \textproc{Burg}((y_n)_{1 \leq n \leq N})$\vspace{1.5 mm}
			\If{$\sum\limits_{m=1}^{d-1} (d-m) \text{atanh}^2\left(\left|\frac{\nu_m - \mu_m}{1-\overline{\mu_m} \nu_m}\right|\right) \leq \epsilon$}\vspace{1.5 mm}
			\State \textbf{break}
			\Else
			\State $\mu \gets \nu$
			\EndIf
			\EndFor
			\State \Return $R_{-1}$
			\EndFunction
		\end{algorithmic}
	\end{samepage}
\end{onehalfspace}
\vspace{2 mm}

Here $\textproc{Burg}((x_n)_{1\leq n\leq N})$ designates the multisegment Burg procedure, returning the residual error power $\sigma$ and the $d-1$ Schur coefficients $(\mu_m)_{1 \leq m \leq d-1}$ \cite{haykin1982maximum}; $\textproc{Trench}(\sigma^2,(\mu_m)_{1 \leq d-1})$ designates the Trench algorithm, returning the inverse covariance matrix $\Sigma^{-1}$ corresponding to the residual error $\sigma$ and the Schur coefficients $(\mu_m)_{1 \leq d-1}$.\\
Note that there is no theoretical guaranty of convergence of the Burg-Tyler algorithm; however numerical convergence has been observed so far on every experiment performed by the authors, on both simulated and real data.

\subsection{Maximum likelihood under a unique scaled radial distribution}
\label{ssec:scaled_radial}

We shall now study the specific case in which $Q$ is modeled by $dQ(r) = \sigma dQ_0(\sigma r)$, $Q_0$ being a fixed probability distribution absolutely continuous with respect to $\mu_{\mathbb{R}_+}$, with density $q$.\\
Putting $\Sigma = \sigma R$, the likelihood can be expressed up to constant terms as:

\begin{equation}
\label{eq:logl_scaled_radial}
\begin{split}
l(Q,R|\delta_{[x_0]}^{\mathbb{K}^d}) = &
-\frac{c_\mathbb{K}}{2}\log\left|\Sigma\right|\\
&  -H\left(\delta_{\sqrt{{x_0}^\dagger \Sigma^{-1} x_0}}^{\mathbb{R}_+}|Q_0\right)\\
&-\frac{c_\mathbb{K}(d-1)}{2} \log\left({x_0}^\dagger \Sigma^{-1} x_0\right)
\end{split}
\end{equation}

This can be expressed, up to a constant, as:
\begin{equation}
\label{eq:logl_scaled_radial_elliptical}
\begin{split}
l(Q,R|\delta_{[x_0]}^{\mathbb{K}^d}) = &
-\frac{c_\mathbb{K}}{2}\log\left|\Sigma\right|\\
&  +\frac{1}{2} \log\left(q\left({x_0}^\dagger \Sigma^{-1}x_0\right)\right)\\
&-\frac{c_\mathbb{K}(d-1)}{2} \log\left({x_0}^\dagger \Sigma^{-1} x_0\right)
\end{split}
\end{equation}

which we shall simplify as:

\begin{equation}
\label{eq:logl_scaled_radial_g}
l\left(Q,R|\delta_{[x_0]}^{\mathbb{K}^d}\right) = 
-\frac{c_\mathbb{K}}{2}\log\left|\Sigma\right|+\frac{1}{2} g\left({x_0}^\dagger \Sigma^{-1}x_0\right)
\end{equation}

with $g$ being defined by:

\begin{equation}
\label{eq:g_def}
g(t) = \frac{1}{2}(c_\mathbb{K}(d-1)\log(t)-\log(q(t)))
\end{equation}

Assuming that $g$ is $\mathcal{C}^1$, this expression can be differentiated and summed over several samples $\left(x_n\right)_{1 \leq n \leq N}$, which produces the following maximum likelihood equation:

\begin{equation}
\label{eq:ml_scaled_radial}
\text{tr}\left(\left(\Sigma-\frac{1}{N}\sum\limits_{n=1}^N g'\left({x_n}^\dagger \Sigma^{-1}x_n\right)x_n {x_n}^\dagger\right)d\left(\Sigma^{-1}\right)\right)		
\end{equation}

Under the condition that $g'\geq 0$  and verifies some necessary conditions detailed in \cite{ollila2003robust} and \cite{marona1998robust}, one can use the following standard M-estimator of the covariance matrix:

\begin{onehalfspace}
\begin{samepage}
\begin{algorithmic}[1]
	\Function{M\_cov}{$g',\left(x_n\right)_{1 \leq n \leq N},\epsilon,K_{max}$}
		\State $\Sigma \gets I$
		\Repeat
			\State $\Sigma_{-1} \gets \Sigma^{-1}$\vspace{1.5 mm}
			\State $\Sigma \gets \frac{1}{N}\sum\limits_{n=1}^N g'\left({x_n}^\dagger {\Sigma_{-1}} x_n\right)x_n{x_n}^\dagger$\vspace{1.5 mm}
		\Until{$\text{tr}\left(\left(\Sigma_{-1} \Sigma - I\right)^2\right) \leq \epsilon$}\vspace{1.5 mm}
		\State \Return $\Sigma_{-1}$
	\EndFunction
\end{algorithmic}
\end{samepage}
\end{onehalfspace}
\vspace{2 mm}

For constrained optimizations, one can use the M version of a known estimator for gaussian models:\\

\begin{onehalfspace}
\begin{samepage}
\begin{algorithmic}[1]
	\Function{M\_of}{$g'$,$e$,${x_n}_{1 \leq n \leq N},\epsilon,K_\text{max}$}
		\State $\Sigma_{-1} \gets I$
		\For{$k$ \textbf{from} $1$ \textbf{to} $K_\text{max}$}
			\For{$n$ \textbf{from} $1$ \textbf{to} $N$}
				\State $y_n \gets \sqrt{g'\left({x_n}^\dagger \Sigma_{-1} x_n\right)}x_n$
			\EndFor\vspace{1.5 mm}
			\State $\Sigma \gets e\left(\left(y_n\right)_{1 \leq n \leq N}\right)$\vspace{1.5 mm}
			\If{$\text{tr}\left({\left(\Sigma_{-1}\Sigma-I\right)}^{2}\right) \leq \epsilon$}\vspace{1.5 mm}
				\State \textbf{break}
			\Else
				\State $\Sigma_{-1} \gets \Sigma^{-1}$
			\EndIf
		\EndFor
		\State \Return $\Sigma_{-1}$
	\EndFunction
\end{algorithmic}
\end{samepage}
\end{onehalfspace}
\vspace{2 mm}

Now if one does not have the condition $g' \geq 0$, it is necessary to resort to geodesic shooting on the manifold of positive definite matrices. This corresponds to the following procedure:

\begin{onehalfspace}
\begin{samepage}
\begin{algorithmic}[1]
	\Function{M\_exp\_cov}{$g',\left(x_n\right)_{1 \leq n \leq N},\epsilon,K_{max}$}
		\State $\Sigma \gets I$
		\Repeat
			\State $\Sigma_{\frac{1}{2}} \gets \sqrt{\Sigma}$
			\State $\Sigma_{-\frac{1}{2}} \gets {\Sigma_{\frac{1}{2}}}^{-1}$\vspace{1.5 mm}
			\State $S \gets \frac{1}{N}\sum\limits_{n=1}^N g'\left({x_n}^\dagger {\Sigma_{-1}} x_n\right)x_n{x_n}^\dagger$\vspace{1.5 mm}
			\State $\Sigma \gets \Sigma_{\frac{1}{2}}\exp\left(\Sigma_{-\frac{1}{2}} S \Sigma_{-\frac{1}{2}}-I\right)\Sigma_{\frac{1}{2}}$\vspace{1.5 mm}	
		\Until{$\text{tr}\left(\left(\Sigma_{-\frac{1}{2}} \Sigma \Sigma_{-\frac{1}{2}} - I\right)^2\right) \leq \epsilon$}\vspace{1.5 mm}	
		\State \Return $\Sigma_{-\frac{1}{2}}$
	\EndFunction
\end{algorithmic}
\end{samepage}
\end{onehalfspace}
\vspace{2 mm}

Here $\sqrt{S}$ designates the unique positive definite square root of the positive definite matrix $S$, and $\exp(A)$ designates the matrix exponential of any matrix $A$.

Unfortunately there is no way of extending this procedure for other known estimators for constrained estimation problems.

Note  that in the case of a gaussian in $\mathbb{R}^d$, one has $g'(t) = 1$ and therefore the corresponding estimator is the sample covariance.\\

In the case of a gaussian in $\mathbb{C}^d$, one has, under the circularity hypothesis: $g'(t) = 1-\frac{1}{2t}$.\\

Thus one can use the geodesic procedure to derive the following algorithm:

\begin{onehalfspace}
\begin{samepage}
\begin{algorithmic}[1]
	\Function{cg\_cov}{$\left(x_n\right)_{1 \leq n \leq N},\epsilon,K_{max}$}
		\State $\Sigma \gets \frac{1}{N}\sum\limits_{n=1}^{N} x_n {x_n}^\dagger$\vspace{1.5 mm}
		\Repeat
			\State $\Sigma_{\frac{1}{2}} \gets \sqrt{\Sigma}$
			\State $\Sigma_{-\frac{1}{2}} \gets {\Sigma_{\frac{1}{2}}}^{-1}$\vspace{1.5 mm}
			\State $S \gets \frac{1}{\left(1-\frac{1}{2d}\right)N}\sum\limits_{n=1}^N \left(1-\frac{1}{2{x_n}^\dagger {\Sigma_{-1}} x_n}\right)x_n{x_n}^\dagger$\vspace{1.5 mm}
			\State $\Sigma \gets \Sigma_{\frac{1}{2}}\exp\left(\Sigma_{-\frac{1}{2}} S \Sigma_{-\frac{1}{2}}-I\right)\Sigma_{\frac{1}{2}}$\vspace{1.5 mm}
		\Until{$\text{tr}\left(\left(\Sigma_{-\frac{1}{2}} \Sigma \Sigma_{-\frac{1}{2}} - I\right)^2\right) \leq \epsilon$}\vspace{1.5 mm}
		\item \textbf{return} $\Sigma_{-\frac{1}{2}}$
	\EndFunction
\end{algorithmic}
\end{samepage}
\end{onehalfspace}
\vspace{2 mm}

\section{Target detection  in a CES background of unknown radial distribution}
\label{sec:detec}

Once the correlation structure of a data set containing a CES distributed background has been estimated, an application of particular interest is the detection of super-imposed signals in such a background. In this article, we shall concentrate on the detection of 1-dimensional deterministic signals of unknown amplitudes.\\

\subsection{Mormalized matched filter}
\label{ssec:nmf}

One performs a test between the two following hypothesis:
\begin{itemize}
	\item $H_0$ : $x = x_R$ : $x$ follows a CES distribution of unknown radial distribution, of correlation matrix $R$
	\item $H_1$ : $x = x_R+\alpha s$ : $x$ follows a CES ditribution of unknown radial distribution, or correlation matrix $R$, centered around a deterministic signal of unknown amplitude $\alpha s$ ($\alpha \in \mathbb{C},s \in \mathbb{C}^d$).
\end{itemize}

By taking the maximum likelihood under the hypothesis $H_1$ and $H_0$ over possible values of $\alpha$ and of the radial distribution $Q$ and by taking their difference, one obtains the following generalized likelihood ratio test:

\begin{equation}
\begin{split}
\text{NMF}(x|s,R) =& \max_{\alpha \in \mathbb{C},Q} l_{H_1}(x|Q,\alpha)\\
 &- \max_{Q} l_{H_0}(x|Q) \geq T(p_{fa})
\end{split}
\end{equation}

Using the expression of the concentrated likelihood for the sampling distribution $S(y) = \delta_x^{\mathbb{C}^d}(y)$:

\begin{equation}
\begin{split}
\text{NMF}(x|s,R) = &\max_{\alpha \in \mathbb{C}} -(d-1)\log\left(\frac{(x-\alpha s)^\dagger R^{-1} (x - \alpha s)}{x^\dagger R^{-1} x}\right) \\
&= -(d-1)\log\left(1-\frac{|s^\dagger R^{-1} x|^2}{x^\dagger R^{-1} x s^\dagger R^{-1} s}\right)
\end{split}
\end{equation}

This is therefore equivalent to:

\begin{equation}
\frac{|s^\dagger R^{-1} x|^2}{x^\dagger R^{-1} x s^\dagger R^{-1} s} \geq T(p_{fa})
\end{equation}

This test has been derived using different methods by several authors, and thus has appeared with different names\cite{scharf1994matched} \cite{scharf1996adaptive} \cite{conte2002recursive} \cite{gini1997sub} \cite{gini2002covariance}; we have chosen the denomination of normalized matched filter in this paper (NMF).

It should be noted that under the $H_0$ hypothesis , the output of the test $\text{NMF}(x|R)$ follows a $\chi^2(2)$ distribution \cite{scharf1994matched}, and is thus CFAR with respect to the dimension, the correlation matrix and the target signal.\\

In practice the correlation matrix $R$ is not known and a prior estimate is used instead. The corresponding test is then called adaptive normalized matched filter (ANMF). It is still CFAR with respect to the correlation matrix and the target signal if one uses Tyler's estimate for prior estimation of the correlation matrix \cite{mahot2013asymptotic} \cite{frontera2013false}; however it is not CFAR with respect to the space dimension anymore.\\

Note that if one chooses to express the likelihood without taking into account the circular symmetry of the sampling distribution, one gets the same expression up to a scale factor:
\begin{flalign*}
\text{NMF}_\phi(x|s,R) = &&
\end{flalign*}
\begin{equation}
 -\left(d-\frac{1}{2}\right)\log\left(1-\frac{|s^\dagger R^{-1} x|^2}{x^\dagger R^{-1} x s^\dagger R^{-1} s}\right)
\end{equation}

\subsection{Multi-channel adaptive normalized matched filter}
\label{ssec:multi_nmf}

The test can be easily extended to a multi-channel scenario. One supposes that there are $K$ uncorrelated signal subspaces of dimensions $(d_k)_{1 \leq k \leq K}$, each of which being spanned by a CES background of known correlation matrix $R_k$ and unknown radial distribution $Q_k$. One wants to detect a signal which can be represented on each subspace by the vector $s_k$, with an unknown amplitude $\alpha_k$. The GLR test in this situation is then given by:

\begin{flalign*}
\text{NMF}(\left(x_k\right)_{1 \leq k \leq K}|\left(s_k\right)_{1 \leq k \leq K},\left(R_k\right)_{1 \leq k \leq K}) &&
\end{flalign*}
\begin{equation}
\begin{split}
=&\sum\limits_{k=1}^K \max_{\alpha_k \in \mathbb{C}} -(d_k-1)\log\left(\frac{\left(x_k-\alpha_k s_k\right)^\dagger {R_k}^{-1} \left(x_k - \alpha_k s_k\right)}{{x_k}^\dagger {R_k}^{-1} x_k}\right)\\
=& -\sum\limits_{k=1}^K (d_k-1)\log\left(1-\frac{|{s_k}^\dagger {R_k}^{-1} x_k|^2}{{x_k}^\dagger {R_k}^{-1} x_k {s_k}^\dagger {R_k}^{-1} s_k}\right) \geq T(p_{fa})
\end{split}
\end{equation}

The output $\text{NMF}\left(\left(x_k\right)_{1 \leq k \leq K}|\left(s_k\right)_{1 \leq k \leq K},\left(R_k\right)_{1 \leq k \leq K}\right)$ of this test follows a $\chi^2(2K)$ distribution under the $H_0$ hypothesis, and is therefore CFAR with respect to the subspace dimensions, the correlation matrices and the target signal on each subspace.\\

In practice, the correlation matrices are not known and prior estimates are used instead. The corresponding test is still CFAR with respect to the correlation matrices and target signals if one uses Tyler's estimates as prior estimates of the correlation matrices; it is not, however, CFAR with respect to subspace dimensions in this case.\\

Again one could also derive the test without taking into account a circular symmetry of the sampling distribution:
\begin{flalign*}
\text{NMF}_\phi(\left(x_k\right)_{1 \leq k \leq K}|\left(s_k\right)_{1 \leq k \leq K},\left(R_k\right)_{1 \leq k \leq K}) =&&
\end{flalign*}
\begin{equation}
 -\sum\limits_{k=1}^K \left(d_k-\frac{1}{2}\right)\log\left(1-\frac{|{s_k}^\dagger {R_k}^{-1} x_k|^2}{{x_k}^\dagger {R_k}^{-1} x_k {s_k}^\dagger {R_k}^{-1} s_k}\right) \geq T(p_{fa})
\end{equation}

However this test shows a slight degradation of performances. Moreover it loses the CFAR property with respect to subspace dimensions.

\subsection{GLR for a scaled fixed radial distribution}
\label{ssec:detec_scaled_radial}

Let us now suppose that the radial distribution is fixed. Using the matrix $\Sigma$ and the function $g$ previously introduced, on has:

\begin{equation}
\begin{split}
\text{GLR}_g(x|s,\Sigma) =& \max_{\alpha \in \mathbb{C}} -g\left({(x-\alpha s)}^\dagger \Sigma^{-1}(x-\alpha s)\right)\\ &+g\left(x^\dagger \Sigma^{-1}x\right)
\end{split}
\end{equation}

If $g$ is $\mathcal{C}^1$, this maximum verifies one of the two following properties:

\begin{equation}
\left\{\begin{array}{lc}
g'\left({(x-\alpha s)}^\dagger \Sigma^{-1} (x-\alpha s)\right) = 0&\mbox{(1)}\\
\alpha = \frac{s^\dagger \Sigma^{-1} x}{s^\dagger \Sigma^{-1} s}&\mbox{(2)}
\end{array}\right.
\end{equation}

The likelihood ratio test of points associated with the first condition depends only on the value of $g'$ at a specific point. Thus the generalized likelihood ratio test in this case is given by:\\

\begin{flalign*}
\text{GLR}_g(x|s,\Sigma) &&
\end{flalign*}
\begin{equation}
\begin{split}
=& g(\tau)\\
&-\min g\left(\left({g'}^{-1}\left(\left\{0\right\}\right) \cap \left[\tau-m;+\infty\right[\right)\cup \left\{\tau-m\right\}\right)
\end{split}
\end{equation}

with $\tau$ and $m$ defined by:

\begin{equation}
\left\{\begin{array}{lc}
\tau = x^\dagger \Sigma^{-1} x& \mbox{(1)}\\
m = \frac{{\left|s^\dagger\Sigma^{-1}x\right|}^2}{s^\dagger \Sigma^{-1}s}& \mbox{(2)}
\end{array}\right.
\end{equation}

This test can be extended in a multi-channel setting by simply summing its output for each channel, as in the case of the $\text{NMF}$ test.\\

Note that if $g'>0$ which is a quite common assumption, the test simplifies to:

\begin{equation}
\text{GLR}_g(x|s,\Sigma) = g(\tau)-g(\tau-m)
\end{equation}

In the particular example of gaussian distributions, this does not coincide with the usual matched filter, which is obtained when the circular symmetry of the sampling distribution is not taken into account.
Indeed the test becomes:
\begin{flalign*}
\text{GLR}_\text{cg}(x|s,\Sigma) &&
\end{flalign*}
\begin{equation}
 = \left\{\begin{array}{lc}
\tau-\frac{1}{2}\log(\tau)-\frac{1}{2}\left(1+\log(2)\right) & \mbox{if } \tau-m\leq\frac{1}{2}\\
m+\frac{1}{2}\log\left(1-\frac{m}{\tau}\right) & \mbox{if } \tau-m>\frac{1}{2}
\end{array}\right.
\end{equation}

with $\tau$ and $m$ defined as previously.\\

This test gives almost identical performances to the matched filter when one knowns the true covariance matrix, although it is slightly better in a multi-channel setting. Moreover in an adaptive context, the combination cg\_cov-$\text{GLR}_\text{cg}$ performs better than the combination of sample covariance and matched filter (SCM-MF).

\section{Simulations}
\label{sec:simu}

We now show some simulation results of detectors in a known covariance background and of adaptive detectors using various estimators introduced in this article.

\subsection{Signal detection in a known covariance background}
\label{ssec:simu_known_background}

We first focus on the case in which the covariance matrix of the background noise is known. The simulation results are shown in a multi-channel scenario, as a function of the SNR of the target signal.\\
The channels are chosen with dimensions 2,4,8 and 16; the noise on each channel is generated as a white gaussian noise of unit variance.\\
The target signal on each channel is generated as a complex centered circular 1-dimensional gaussian signal aligned with the test signal $s_k$, whose variance $\sigma$ on each channel $k$ is such that: 
\begin{displaymath}
10\log_{10}(\sigma_k) = d_k\text{SNR}
\end{displaymath}

The detection thresholds are defined in order to have a false alarm rate of $10^{-4}$.\\

\begin{figure}[H]
	\includegraphics[width=\linewidth]{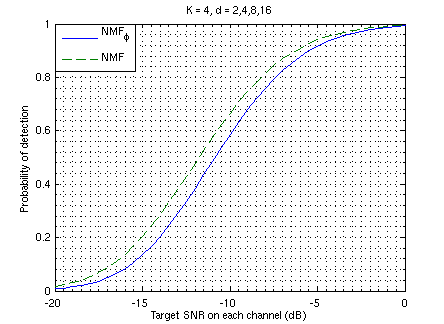}
	\caption{Detection capability of $\text{NMF}_\phi$ and NMF tests in a multi-channel scenario}
	\label{fig_norm_comp}
\end{figure}

\begin{figure}[H]
	\includegraphics[width=\linewidth]{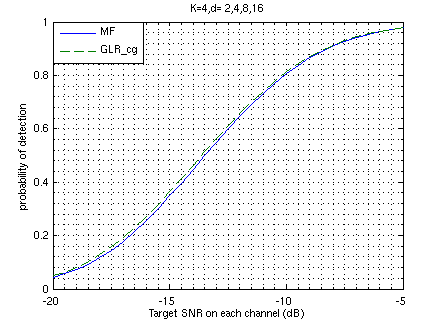}
	\caption{Detection capability of Matched filter (MF) and $\text{GLR}_\text{cg}$ tests in a multi-channel scenario}
	\label{fig_gauss_comp}
\end{figure}

Figure~\ref{fig_norm_comp} shows the probability of detection of the NMF and $\text{NMF}_\phi$ tests as a function of the SNR in a multi-channel scenario; figure~\ref{fig_gauss_comp} shows the detection capabilities of the matched filter (MF) and GLR\_cg tests in the same conditions.\\

As can be seen on figures~\ref{fig_norm_comp} and~\ref{fig_gauss_comp}, taking into account the circular symmetry of signals gives a slight improvement of performances on the tested scenario.

\subsection{Adaptative signal detection}

We now focus on the case in which the covariance matrix of the background noise is unknown. The simulation results are shown in a single channel scenario of dimension $d=8$, as a function of the SiNR of the target signal.\\
The target signal is generated as a complex centered circular 1-dimensional gaussian signal aligned with the test signal $s$, whose variance $\sigma$ is such that:

\begin{displaymath}
10 \log_{10}(\sigma) = \text{SiNR}
\end{displaymath}

The detection thresholds are again defined to have a false alarm rate of $10^{-4}$.\\

Prior covariances are estimated with $N = 22$ independently drawn noise samples.\\

\begin{figure}[H]
	\includegraphics[width=\linewidth]{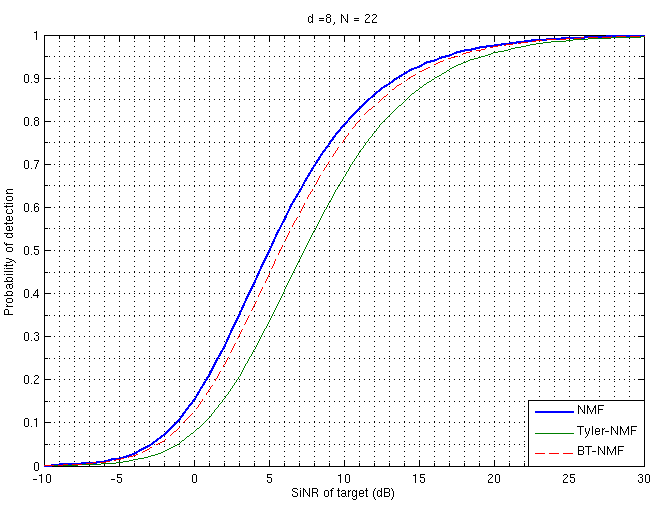}
	\caption{Detection capability of \textproc{Tyler}-NMF and \textproc{BT}-NMF tests}
	\label{fig_anorm_comp}
\end{figure}

\begin{figure}[H]
	\includegraphics[width=\linewidth]{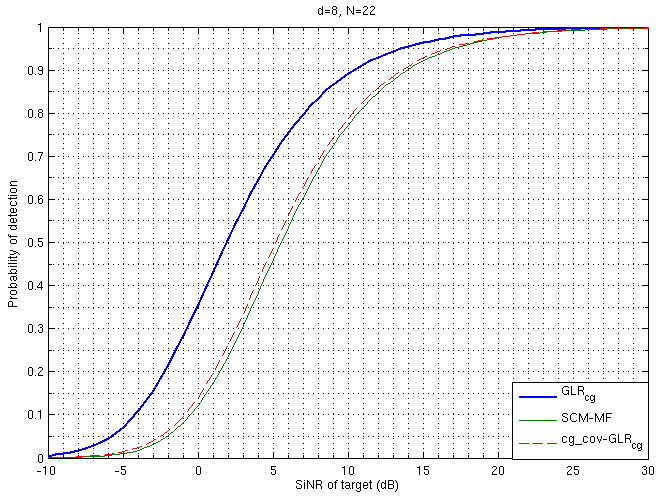}
	\caption{Detection capability of SCM-MF and \textproc{cg\_cov}-$\text{GLR}_\text{cg}$}
	\label{fig_agauss_comp}
\end{figure}

Figure~\ref{fig_anorm_comp} shows the probability of detection of the \textproc{Tyler}-NMF and \textproc{BT}-NMF adaptive detectors, as well as the NMF using the real covariance matrix for comparison. \\
Figure~\ref{fig_agauss_comp} shows the probability of detection of the SCM-MF and \textproc{cg\_cov}-$\text{GLR}\_\text{cg}$ adaptive detectors, as well as the \textproc{cg\_cov} using the real covariance matrix for comparison.

As can be seen on figures~\ref{fig_anorm_comp} and~\ref{fig_agauss_comp}, all adaptive detectors suffer from losses compared to the same detectors when the covariance is known. These losses can be minimized in some cases by making further assumptions on the signal model, such as stationarity (figure~\ref{fig_anorm_comp}) and circularity (figure~\ref{fig_agauss_comp}).

\section{conclusion}

We have presented a complete and coherent theory of maximum likelihood estimation problems under the usual i.i.d hypothesis, which can be equivalently seen as a minimum information loss problem. This can be used to properly derive ML estimators in the case of potentially degenerate distribution models which do not admit a density for the standard Lebesgues measure. This background is used to rederive Tyler's estimators and other types of M-estimators as likelihood maximizers.\\
This formalism also allows to take into account a priori known information about sampling distributions; we show in particular how this can be used to take into account known symmetries of a problem. In particular, we have shown how the circularity of complex signals can be properly accounted for. The corresponding estimators are derived, as well as the statistical tests for signal detection in complex elliptically distributed background and their correct extension to a multichannel setting.



\ifCLASSOPTIONcaptionsoff
  \newpage
\fi



%

\bibliographystyle{plain}
\bibliography{CES_biblio}



%

\begin{IEEEbiographynophoto}{Christophe Culan}
	was born in Villeneuve-St-Georges, France, on November $\text{1}^\text{st}$, 1988. He received jointly the engineering degree from Ecole Centrale de Paris (ECP), France and a master's degree in Physico-informatics from Keio University, Japan, in 2013.\\
	He was a researcher in applied physics in Itoh laboratory from 2013 to 2014, specialized in quantum information and quantum computing, and has contributed to several publications related to these subjects.\\
	He currently holds a position as a research engineer in Thales Air Systems, Limours, France, in the Advanced Radar Concepts division. His current research interests include statistical signal and data processing, robust statistics, machine learning and information geometry.
\end{IEEEbiographynophoto}

\begin{IEEEbiographynophoto}{Claude Adnet}
	was born in Aÿ, France, in 1961. He received the DEA degree in signal processing and Phd degree in 1988 and 1991 respectively, from the Institut National Polytechnique de Grenoble (INPG), Grenoble France.
	Since then, he has been working for THALES Group, where he is now Senior Scientist. His research interests include radar signal  processing and radar data processing.
\end{IEEEbiographynophoto}





\end{document}